\magnification=1200
\input amssym.def
\input epsf

\def \ft {FT\ }
\def \dft {$\Delta$FT\ }
\def \dn {{$\Delta_n$}}
\def\sqr#1#2{{\vcenter{\vbox{\hrule height.#2pt
     \hbox{\vrule width.#2pt height#1pt \kern#1pt
     \vrule width.#2pt} \hrule height.#2pt}}}}
\def\square{\ ${\mathchoice\sqr34\sqr34\sqr{2.1}3\sqr{1.5}3}$}

\def \BZ {{\Bbb Z}}
\def \LL {{\cal L}}
\def \cl#1{\lfloor #1 \rceil}
\def \implies {\Rightarrow}

\centerline {\bf BRAID COMMUTATORS AND DELTA FINITE-TYPE INVARIANTS}
\bigskip
\centerline {Theodore B. Stanford \footnote \dag 
{Research supported in part by the Naval Academy Research Council}}
\centerline {Mathematics Department}
\centerline {United States Naval Academy}
\centerline {572C Holloway Road}
\centerline {Annapolis, MD\ \ 21402}
\medskip
\centerline {\tt stanford@nadn.navy.mil}

\bigskip
\bigskip
\noindent
{\bf ABSTRACT}
\bigskip
Delta finite-type invariants are defined analogously to
finite-type invariants, using delta moves instead of 
crossing changes.  We show that they are closely related
to the lower central series of the commutator subgroup
of the pure braid group.
\bigskip
\bigskip

\bigskip
\noindent
{\bf 0. INTRODUCTION}
\bigskip

We consider in this paper
delta finite-type
invariants of knots and links (\dft invariants).
In the case of links, these are the same invariants
as defined by Mellor~[4].
We shall prove some properties of
\dft invariants which closely resemble properties
of finite-type invariants in the usual sense
(\ft invariants).  In the same
way that \ft invariants are based on sets of crossing changes
in knot or link diagrams, \dft invariants are based on sets of
delta moves in knot or link diagrams.  In the same way that
\ft invariants are closely related to $\gamma_n (P)$, the
the lower central series of the pure braid group $P$, \dft
invariants are closely related to $\gamma_n (P^\prime)$, the
lower central series of the commutator
subgroup of $P$. Here $P = P_k$, the pure braid group on
$k$ strands.

We will show that two knots have matching
\dft invariants of order $<n$ if and only if they are
equivalent modulo $\gamma_n (P^\prime)$, just as two knots
have matching \ft invariants of order $<n$ if and only if
they are equivalent modulo $\gamma_n (P)$ 
(see [8]).  Since
$\gamma_n (P^\prime) \subset \gamma_{2n} (P)$, it follows
that an \ft invariant of order $<2n$ is a \dft invariant of
order $<n$.  It may turn out that all \dft invariants of
order $<n$ occur this way.  However, there is a difference
between $\gamma_n (P)$ and $\gamma_n (P^\prime)$ which 
may make the
\dft invariants more than just a relabeling of the
\ft invariants:	\  For any $k$ and any $n$, the
quotient group $\gamma_n (P_k)/\gamma_{n+1} (P_k)$ is
finitely-generated, whereas even $\gamma_1
(P_3^\prime)/\gamma_2 (P_3^\prime)$ is not finitely-generated.

For links of more than one component, there are known
advantages to working with \dft invariants instead of
\ft invariants.   Murakami and 
Nakanishi~[5] showed 
that two links are equivalent by a sequence
of delta moves if and only if they have the same pairwise
linking numbers.  Thus, \dft invariants detect linking number
right at order $0$, and so essentially each link homology
class gets its own set of invariants of higher orders.
Invariants which are not \ft over all links may be \dft within
one homology class. For example, 
Mellor~[4] has shown that
Milnor's triple linking number is \dft but not \ft (it is not
even well-defined over all links).  See also 
Appleboim~[1].

\bigskip
\noindent
{\bf Appreciation:} \ 
I would like to thank the organizers of the
``Knots in Hellas'' conference, where this work
was presented.

\bigskip
\noindent
{\bf 1. DEFINITIONS}
\bigskip

The basic terminology here is motivated by the
standard ideas from \ft invariants.  See for example 
Birman~[2], 
Gusarov~[3], 
Ohyama~[6],
or Taniyama~[9].

\medskip
\noindent
{\bf Notation 1.1.} \  
Let $L$ be a link diagram with
$n$ disjoint sets of disjoint disks $S_1,
S_2, \dots S_n$.
Inside each disk, the diagram is supposed
to look like one side
of Figure~1.2.
For any subset $T \subset \{1,2,3, \dots n\}$,
denote by $L_T$ the link obtained by applying
a delta move in each disk in $S_i$ for all
$i \in T$.

\bigskip
\epsfysize = 3truecm
\centerline {\epsffile {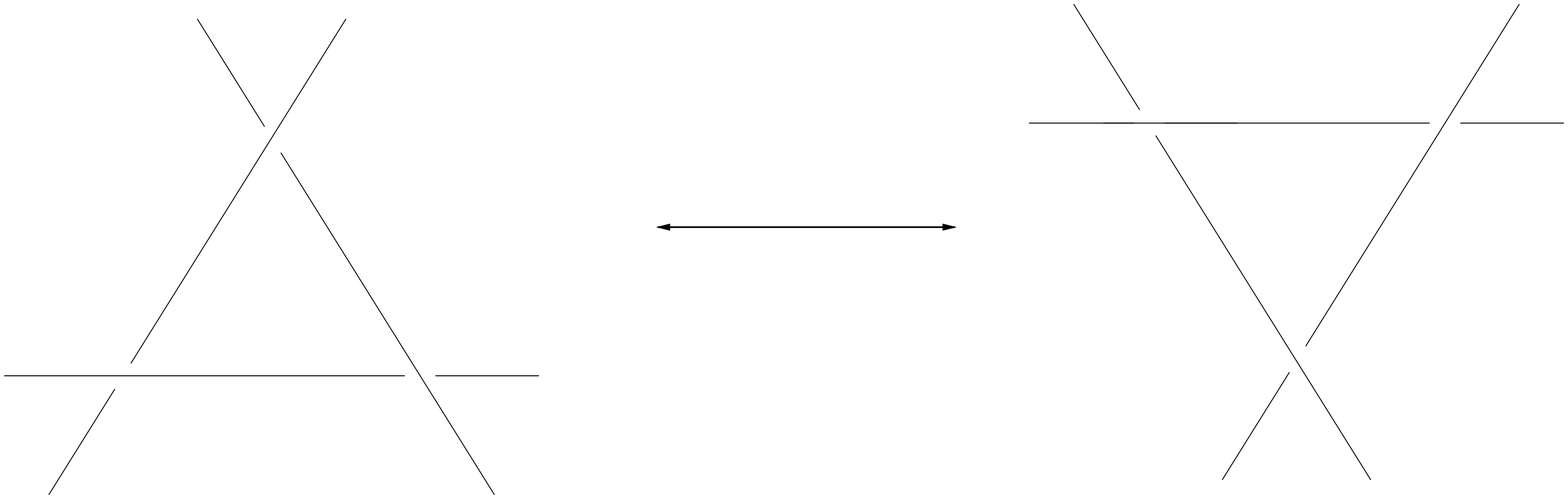}}
\medskip
\centerline 
{Figure 1.2}
\bigskip

\medskip
\noindent
{\bf Definition 1.3.} \ 
Two links $L$ and $L^\prime$
are said to be {\it \dn-equivalent}
if there exists a diagram as
above such that $L_\emptyset = L$ and
$L_T = L^\prime$
for all $T \ne
\emptyset$.  A link \dn-equivalent to
an unlink is said to be {\it \dn-trivial}.

\medskip
I do not know if this is an equivalence
relation in general. If not, then one
can consider the equivalence relation
that it generates.
It will follow from Theorem~2.1 that
the definition does give an
equivalence relation when restricted
to knots.

\medskip
\noindent
{\bf Notation 1.4.} \
Let $\BZ \LL$ be the $\BZ$-module
freely generated by all link types, and
let $C_n$ be the submodule generated
by all linear combinations
$$\sum_{T \subset \{1,2, \dots n\}} (-1)^{|T|} L_T 
\leqno {(1.5)}$$
for all diagrams $L$ as above.

\medskip
\noindent
{\bf Definition 1.6} \ 
Let $v$ be a link invariant taking values
in some abelian group $A$.  Extend $v$ linearly
to $\BZ \LL$.  Then $v$ is said
to be {\it \dft of order $< n$} if it
vanishes on $C_n$.

\vfil
\eject

\noindent
{\bf Proposition 1.7.} \
{\it If two links $L$ and $L^\prime$
are \dn-equivalent
and $v$ is a \dft invariant of order $<n$, then
$v (L) = v(L^\prime)$. }

\medskip
\noindent
{\it Proof:} \   
Choose a diagram for $L$ as in Definition~1.3.
Use this diagram to form the alternating sum as in 
Notation~1.4.
\square

\medskip
\noindent
{\bf Notation 1.8.} \
If $b \in B_k$ is a braid, then we denote by $\cl b$
the usual closure of $b$, connecting the top ends to
the bottom ends to make a knot or link.  This extends to
a $\BZ$-linear map from $\BZ B_k$ to $\BZ \LL$.  Set
$t_k = \sigma_{k-1}^{-1} \sigma_{k-2}^{-1} \dots \sigma_1^{-1}$.  

\medskip
If $p \in P_k$, then $\cl {p t_k}$ is a knot,
and any knot may be written this way.
Also, $t^{-1}_k p t_k$ is $p$ ``shifted'' to
the right. These are very useful facts.

\medskip
\noindent
{\bf Definition 1.9.} \ 
Two links $L$ and $L^\prime$
are said to be
{$\gamma_n (P^\prime)$-equivalent} if there
exists a positive integer $k$ and braids
$p, b \in B_k$ with $L = \cl b$,
$L^\prime = \cl {pb}$, and
$p \in \gamma_n (P_k^\prime)$.  Here
$G^\prime = [G, G]$ and 
$\gamma_n (G) = [G, \gamma_{n-1} (G)]$,
with $\gamma_1 (G) = G$ for a group $G$.

\medskip
The groups $\gamma_n (P_k^\prime)$ form a subcoherent
sequence of pure braid subgroups, and therefore 
Definition~1.9 gives an equivalence
relation.  See~[8].

\medskip
\noindent
{\bf Definition 1.10.} \ 
A knot invariant $v$ is said to be {\it additive} if
$v (K \# K^\prime) = v(K) + v(K^\prime)$.

\bigskip
\noindent
{\bf 2. THEOREMS}
\bigskip

Just as with \ft invariants, there are several different
ways of characterizing the knots whose 
whose \dft invariants match up to a given finite order.

\bigskip
\noindent
{\bf Theorem 2.1.} \ 
{\it For two knots $K$ and $K^\prime$,
the following are equivalent:
\medskip

\smallskip
\item {A.} $K$ and $K^\prime$ are 
$\gamma_n (P^\prime)$-equivalent.

\smallskip
\item {B.} $K$ and $K^\prime$ are \dn-equivalent.

\smallskip
\item {C.} For any \dft invariant $v$ of order
$<n$, $v(K) = v(K^\prime)$.

\smallskip
\item {D.} For any additive \dft invariant $v$ of order
$<n$, $v(K) = v(K^\prime)$.}

\medskip
Gusarov~[3] 
showed that the equivalence classes of knots 
with matching \ft invariants of order $<n$ form a group
under the operation of connected sum.  The same is true
for \dft invariants.

\bigskip
\noindent
{\bf Theorem 2.2.} \ 
{\it For all n, the \dn-equivalence
classes of knots form a group under 
connected sum.}

\medskip
Theorems~2.1 and~2.2 apply only to knots. The next
theorem applies to both knots and links.

\bigskip
\noindent
{\bf Theorem 2.3.} \ 
{\it For any link $L$ and any positive
integer $n$, there exist an infinite number
of prime, nonsplit, alternating links 
$L^\prime$ such that $v (L) = v (L^\prime)$
for all \dft invariants of order $<n$.}

\vfil
\eject

\noindent
{\bf 3. PROOFS}
\bigskip

We showed in [8] 
that the $P^{(n)}$-equivalence classes
of knots (where $P^{(n)}$ is the derived series)
form a group under connected sum, so Theorem~2.2
follows from Theorem~2.1 because
$P^{(n+1)} \subset \gamma_n (P^\prime)$.
We showed in 
[7]
how to modify a link $L$ to create an infinite family
of prime, nonsplit, alternating links in the same
$\gamma_n(P)$-equivalence class as $L$, and we showed
how to make this procedure work for $P^{(n)}$-equivalence
of knots in [8].  
The same procedure is easily modified to give a proof
of Theorem~2.3, 
given Propositions~1.7 
and~3.1.

\medskip
\noindent
{\bf Proposition 3.1.} \ 
{\it If two links
are $\gamma_n (P^\prime)$-equivalent, then
they are \dn-equivalent.}

\medskip
\noindent
{\it Proof:} \  
It suffices to show that a braid
in $\gamma_n (P^\prime)$ is \dn-trivial, in
the sense that there exists a braid diagram
$K$ such that Notation~1.1 and 
Definition~1.3
apply.  The proof is by induction on $n$.
We first need to see that any braid in $P^\prime$
is $\Delta_1$-trivial. 
Let $p_{i,j}$ be the standard generator
of $P_k$ which links the $i$th and the $j$th
strands.  Any commutator of form
$[p_{h,i},p_{i,j}] =
p_{h,i} p_{i,j} p_{h,i}^{-1} p_{i,j}^{-1}$.
may be transformed into a trivial braid by
a delta move. Take any standard set of relations
among the $p_{i,j}$, and add in all relations
of the form $[p_{h,i},p_{i,j}]$.
The result is a presentation of $P_k$ abelianized,
and therefore any braid in $P_k^\prime$ may
be undone by delta moves.

For the induction step, consider
a braid $[p,q]$, where 
$p \in \gamma_n (P^\prime)$ and
$q \in P^\prime$.  Write down the obvious
diagram for $[p,q]$, in terms of diagrams
for $p$ and $q$ which have $n$ and $1$
disjoint sets of disks, respectively,
as in Notation~1.1.  The first $n$ sets
of disks in $[p,q]$ will come from the
disks in $p$ together with the disks
in $p^{-1}$. The last set will come from
the disks in $q$ and $q^{-1}$.
\square

\bigskip
Proposition~3.1 gives us
(A) $\implies$ (B) in Theorem~2.1,
Proposition~1.7 gives us
(B) $\implies$ (C), and (C) $\implies$ (D) is
trivial.  We need to show then that (D) $\implies$ (A). 

If the module $C_n \subset \BZ \LL$ is pulled back to
the group ring $\BZ B_k$ via the closure map, the result is
an easily-described ideal:

\bigskip
\noindent
{\bf Proposition 3.2.} \ 
{\it For each positive integer $k$, let $I = I_k \subset \BZ B_k$
be the ideal generated by all $p-1$ with $p \in P_k^\prime$.  
Then the module $C_n$ is generated by all $\cl {x}$ with
$x \in I_k$ for some $k$.}

\medskip
\noindent
{\it Proof:} \  
Let $L$ be a diagram as in Notation~1.1.
Using a standard trick 
from \ft invariants (see Gusarov~[3]), 
we may restrict ourselves to the
relations~(1.5) where each $S_i$
contains a single disk.
We would like to apply an Alexander-type
theorem as in Birman~[2].
Whether or not a delta move can be fit locally into
a braid diagram depends on the pattern of orientations of
the strands.  (If our links aren't oriented, then choose
an arbitrary orienation for each component.)
Each dot in Figure
Figure~3.3 indicates a possible delta move,
and all four delta moves shown have the same effect.
The moves on the left and right ends of the figure
cannot occur locally in a braid because one of the three
strands will be forced to go ``backwards''.  However,
the figure demonstrates that we can instead always choose
moves like the two in the middle, which can occur locally
in a braid.  We may then write $K$ as a
closed braid diagram, and we may write each delta move as
$[p_{h,i},p_{i,j}] - 1$.  We then 
find that the relation determined by $L$ is in fact the
closure of an element of $I^n$.
\square

\bigskip
\epsfysize = 2truecm
\centerline 
{\epsffile {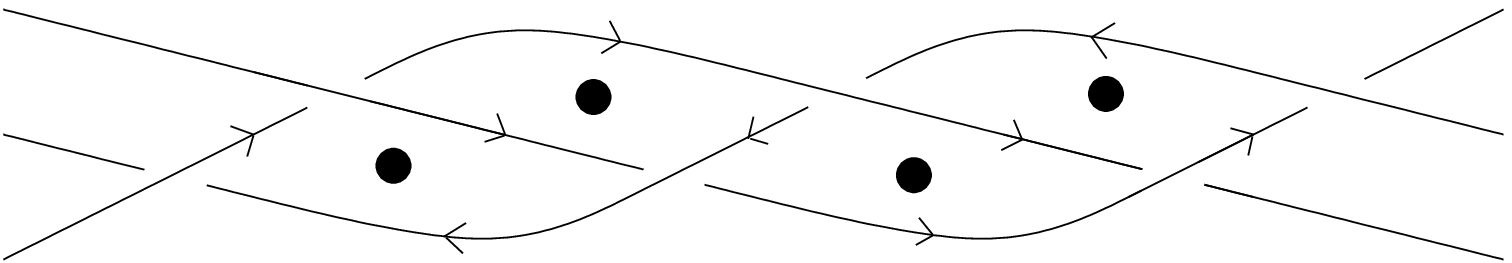}}
\medskip
\centerline
{Figure 3.3}
\bigskip

Now, restricting ourselves to the case of knots, we
may conjugate if necessary and 
see that $C_n$ is generated by expresssions of the form 
$$\cl {(x_1 - 1) (x_2 - 1) \dots (x_n - 1) y t_k}
\leqno {(3.4)}$$
where $x_i \in P_k^\prime$ and $y \in P_k$.
The only thing that prevents us from applying
exactly the same analysis as in Section~2 
of~[8]
is that $y$ is not necessarily in $P_k^\prime$.  This is
fixed as follows.

The expression~(3.4)  
is a $\BZ$-sum of $2^n$ different
knots.  If we take the connected sum of each of these
with the knot $\cl {y^{-1} t_k}$, then we may write
the resulting $\BZ$-sum of knots as 
$$\cl {(x_1 - 1) (x_2 - 1) \dots (x_n - 1)
y (t_{2k}^{-k}y^{-1}t_{2k}^{k}) t_{2k}}
\leqno {(3.5)}$$  If $v$ is an additive
invariant, then it vanishes on~(3.4)  
if and only if it
vanishes on~(3.5).  
Now we use the same idea as 
in~[8],
to slide the $y^{-1}$ around and around until it cancels
with the $y$.  That is, we consider  relations of the 
form 
$$\cl {(x_1 - 1) (x_2 - 1) \dots (x_n - 1)
z y (t_{2k}^{-m}y^{-1}t_{2k}^{m}) t_{2k}}
\leqno {(3.6)}$$
where $0 \le m \le k$, $y \in P_k$, and $z \in P_{2k}^\prime$.
Relation~(3.5) is of this form when $m=k$.
When $m=0$, then $y$ and $y^{-1}$ cancel and we have
the form we need.  Conjugating by 
$t_{2k}^{-m} y t_{2k}^m$, we obtain
$$\cl {(x_1^\prime - 1) (x_2^\prime - 1) \dots (x_n^\prime - 1)
z^\prime y (t_{2k}^{-m+1}y^{-1}t_{2k}^{m-1}) t_{2k}}
\leqno {(3.7)}$$
where 
$x_i^\prime = t^{-m}_{2k} y t_{2k}^m 
x_i t_{2k}^{-m}y^{-1} t_{2k}^m$ and 
$z^\prime = (t^{-m}_{2k} y t_{2k}^m 
z t_{2k}^{-m} y^{-1} t_{2k}^m)
(t^{-m}_{2k} y t_{2k}^m y 
t_{2k}^{-m}y^{-1} t_{2k}^m y^{-1})$.  This is the same
relation as $(3.6)$, and it is written
in the same form except that $m$ has been decreased by one.
Hence a relation of the form (3.4) may
be replaced by a relation of the same form, with $y \in P^\prime$.

\vfil
\eject

\noindent
{\bf REFERENCES}
\medskip

\smallskip
\item {[1]}
E. Appleboim.
{\it Finite type invariants of links with a fixed
linking matrix.}
\hfil \break
{\tt math.GT/9906138}

\smallskip 
\item {[2]} 
J. S. Birman. 
{\it New points of view in knot theory.}
Bulletin of the American Mathematical Society
28 (1993) 253--287.

\smallskip 
\item {[3]} 
M. Gusarov
{\it On $n$-equivalence of knots and invariants of finite degree}.
``Topology of Manifolds and Varieties'', 173--192. 
Advances in Soviet Mathematics 18, 
American Mathematical Society, 1994.

\smallskip
\item {[4]}
B. Mellor.
{\it Finite-type link homotopy invariants II:
Milnor's $\bar {\mu}$ invariants.}
\hfil \break
{\tt math.GT/9812119}.

\smallskip
\item {[5]}
H. Murakami and Y. Nakanishi.
{\it On a certain move generating link-homology}.
Mathematische Annalen 284 (1989),  75--90.

\smallskip
\item {[6]}
Y. Ohyama.
{\it Vassiliev invariants and similarity of knots.}
Proceedings of the AMS 123 (1995), 287--291.

\smallskip 
\item {[7]} 
T. Stanford. 
{\it Braid commutators and Vassiliev invariants}.
Pacific Journal of Mathematics 174 (1996), 269--276.

\smallskip
\item {[8]}
T. Stanford.
{\it Vassiliev invariants and knots modulo pure braid subgroups}.
\hfil \break
{\tt math.GT/9805092}.  

\smallskip 
\item {[9]}
K. Taniyama.
{\it On similarity of knots.}
Gakujutsu Kenkyu, School of Education, Waseda University.
Series of Mathematics 41 (1993), 33--36.

\end